\documentclass{svk_long_en}
\usepackage{epstopdf}
\DeclareGraphicsExtensions{.eps}
\usepackage{hyperref}
\hypersetup{
    colorlinks=false,
    pdfborder={0 0 0},
}

\renewenvironment{abstract}{%
\hfill\begin{minipage}{0.95\textwidth}
\rule{\textwidth}{1pt}}
{\par\noindent\rule{\textwidth}{1pt}\end{minipage}}

\begin{document}
\twocolumn[
\title{Fractal Dimension of Self-Affine Signals: \\Four Methods of Estimation}

\author{Hana Krakovsk\'a\inst{1}, Anna Krakovsk\'a\inst{2}
}

\institute{
FMFI UK,
Mlynská Dolina
842~48~Bratislava
\and
Institute of Measurement Science,
Slovak Academy of Sciences,
Dúbravská cesta 9,
841~04~Bratislava \\ \textit{krakovska@savba.sk}}

\maketitle
\begin{@twocolumnfalse}
\begin{abstract}

This paper serves as a complementary material to a poster presented at the  XXXVI Dynamics Days Europe in Corfu, Greece, on June 6th-10th in 2016.

In this study, fractal dimension ($D$) of two types of self-affine signals were estimated with help of four methods of fractal complexity analysis.

The methods include the Higuchi method for the fractal dimension computation, the estimation of the spectral decay ($\beta$), the generalized Hurst exponent ($H$), and the detrended fluctuation analysis. For self-affine processes, the next relation between the fractal dimension, Hurst exponent, and spectral decay is valid: $D=2-H=\frac{5-\beta}{2}$. Therefore, the fractal dimension can be get from any of the listed characteristics. 

The goal of the study is to find out which of the four methods is the most reliable. For this purpose, two types of test data with exactly given fractal dimensions ($D = 1.2, 1.4, 1.5, 1.6, 1.8$) were generated: the graph of the self-affine Weierstrass function and the statistically self-affine fractional Brownian motion. The four methods were tested on the both types of time series. Effect of noise added to data and effect of the length of the data were also investigated.

The most biased results were obtained by the spectral method. The Higuchi method and the generalized Hurst exponent were the most successful.

\end{abstract}
\end{@twocolumnfalse} \vspace{0.5cm}]

\section{Introduction}

In the last decades the fractal geometry approach of describing the natural phenomena around us has become very popular. It often can better describe and investigate the complex structures in nature, than the traditional Euclidean geometry.  Time series that are of fractal nature have been found in many branches of science, including biology and medicine (e.g. waveforms of EEG or ECG), economy and finance (stock market indices, foreign exchange rates, commodity prices, trading volumes, interest rates...), geology and other physical areas. 
Basically, a fractal is an object that exhibits a repeating pattern across different scales. If the replication is the same at every scale, the fractal is called self-similar. If the pieces repeat themselves only when the axes are magnified by different factors, the fractal is called self-affine.

In this study, we are interested in fractal signals that are self-affine, or self-affine in a statistical sense, which in the latter case means that a rescaled version of a small part of the signal has the same statistical distribution as the larger part. 

More specifically, we will deal with estimation of the so called fractal dimension. Fractal dimension determines the irregularity of the signal. It tells us how ''smooth'' or ''rough'' the curve of the graph is. $D$ can take values between $1$ and $2$. The more the graph fills the plane the closer it approaches the value of $2$. The fractal dimension is a useful characteristic  widely applied in numerous fields. 

When we have self-affine time series, some of the properties as $D$, decay of autocorrelation or power spectrum, persistence, etc. are elegantly connected. 
Above all, the next relation between the fractal dimension, Hurst exponent, and spectral decay holds: $D=2-H=\frac{5-\beta}{2}$ \cite{mande1983}, \cite{eke2002}. As a result, there are several ways to get to the value of the fractal dimension. For instance, instead of directly estimating the dimension, you can derive $D$ from the long-term memory dependence or persistence of time series expressed by Hurst exponent or you can get $D$ from the decay of the spectral power of the function. 

In accordance with the above equations, we investigated the following direct and indirect methods to estimate the fractal dimension:
\begin{enumerate}
\item Higuchi’s method for the computation of $D$
\item Generalized Hurst exponent 
\item Detrended fluctuation analysis (DFA)
\item Spectral decay coefficient ($\beta$)
\end{enumerate}
We chose methods that have succeeded in some earlier comparative studies  \cite{estel2001} \cite{schep2002}.
The methods were tested on the next two types of computer-generated fractal signals with known values of $D$: 
\begin{itemize}
\item Weierstrass function
\item Fractional Brownian motion 
\end{itemize}

This paper is organized as follows.
Section 2 describes the data used for testing.
In Section 3, the four methods of fractal complexity estimation are introduced.
Then the methods are tested to find the most efficient one with regard to the type of data, amount of noise contained in signal, and the length of the data set.
In Section 4 and 5, the results are presented and the findings are summarized.

\section{Data}

\subsection*{Weierstrass function}

Weierstrass function (Wf) has simple notation: 
$$ W(t)=\sum_{n=0}^{\infty} \frac{\cos2^n t}{2^{(2-D)n}}$$
where $0<t<2\pi$ . The parameter $D$ is believed to be exactly the Hausdorff dimension of the resulting curve for $1<D<2$, although this was not rigorously proven yet \cite{falco2004}. 

The Weierstrass function is continuous on the whole domain, though non-differentiable anywhere. It was the first function published in $19^{th}$ century with such unintuitive property. The graph of the function has infinite length and does not smooth out after arbitrary number of magnifications. It is a self-similar fractal, which means, that on a large scale it looks the same as do some smaller parts of it after magnification. 

A graph of the function with $D=1.2$ is considered as a very good approximation of the horizon in the mountains. Its fractal quality also explains why it is difficult to judge how close you are to mountain, since the smaller but closer mountains look the same as bigger ones far away. 

We set different values of $D$ between $1$ and $2$: $1.2, 1.4, 1.5, 1.6, 1.8$. Summation was done for $n$ from 0 to 1021. For each value of $D$, 314160 data points were generated from 0 to 2$\pi.$ 

\begin{figure}
\includegraphics[trim=50 40 50 50,clip,width=0.5\textwidth]{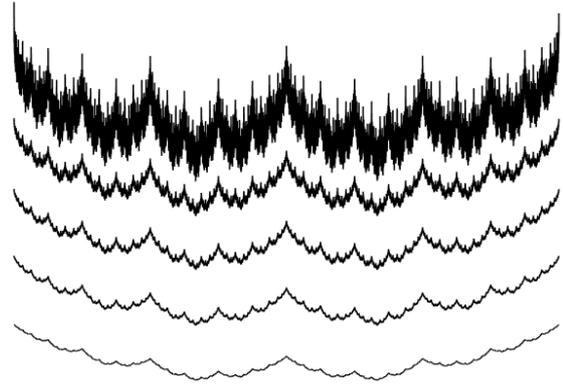}
\caption{Graphs of the Wf with $D=1.2$, $D=1.4$, $D=1.5$, $D=1.6$, $D=1.8$ from bottom to top.} 

\end{figure}  

\subsection*{Fractional Brownian Motion}

Brownian motion (Bm) is a random motion of particles in fluid, which is caused by interactions with atoms and molecules in fluid. As a mathematical model, it is used in different areas of science, and also for a description of price movements in stock market. Fractional Brownian motion (fBm) is a generalization of Bm. It is a non-stationary signal with stationary increments. fBm is a continuous-time Gaussian process with properties depending on Hurst exponent $0<H<1$. While the ordinary Bm has the Hurst exponent equal to $0.5$, for fBm $H$ can vary between $0$ and $1$. Unlike Bm, the increments of fBm are not independent.  For $H>0.5$ the increments of the process are positively correlated, while for $H<0.5$ they are negatively correlated. 
The fractional Brownian motion is statistically self-similar in that its segments are equal in distribution to a longer segments when the latter are properly  rescaled. The variance of the increments is given by $$Var(fBm(t)-fBm(s)) = v*|t-s|^{2H}$$
where $v$ is a positive constant. Fractal dimension can be derived from the relation $D=2-H$. So when generating fBm, we can set the value of the Hurst exponent and thereby also $D$  and generate data with known $D$.  For each of the chosen values of $D$ $(1.2, 1.4, 1.5, 1.6, 1.8)$ we generated $320 000$ data points. To do so we used the Matlab function $wfbm(H,L)$ where we set the Hurst parameter $H$ ($0 < H < 1$) and the length of data $L$. The function is based on algorithm of Abry and Sellan  \cite{abry1996}. 

\begin{figure}
\includegraphics[width=0.5\textwidth]{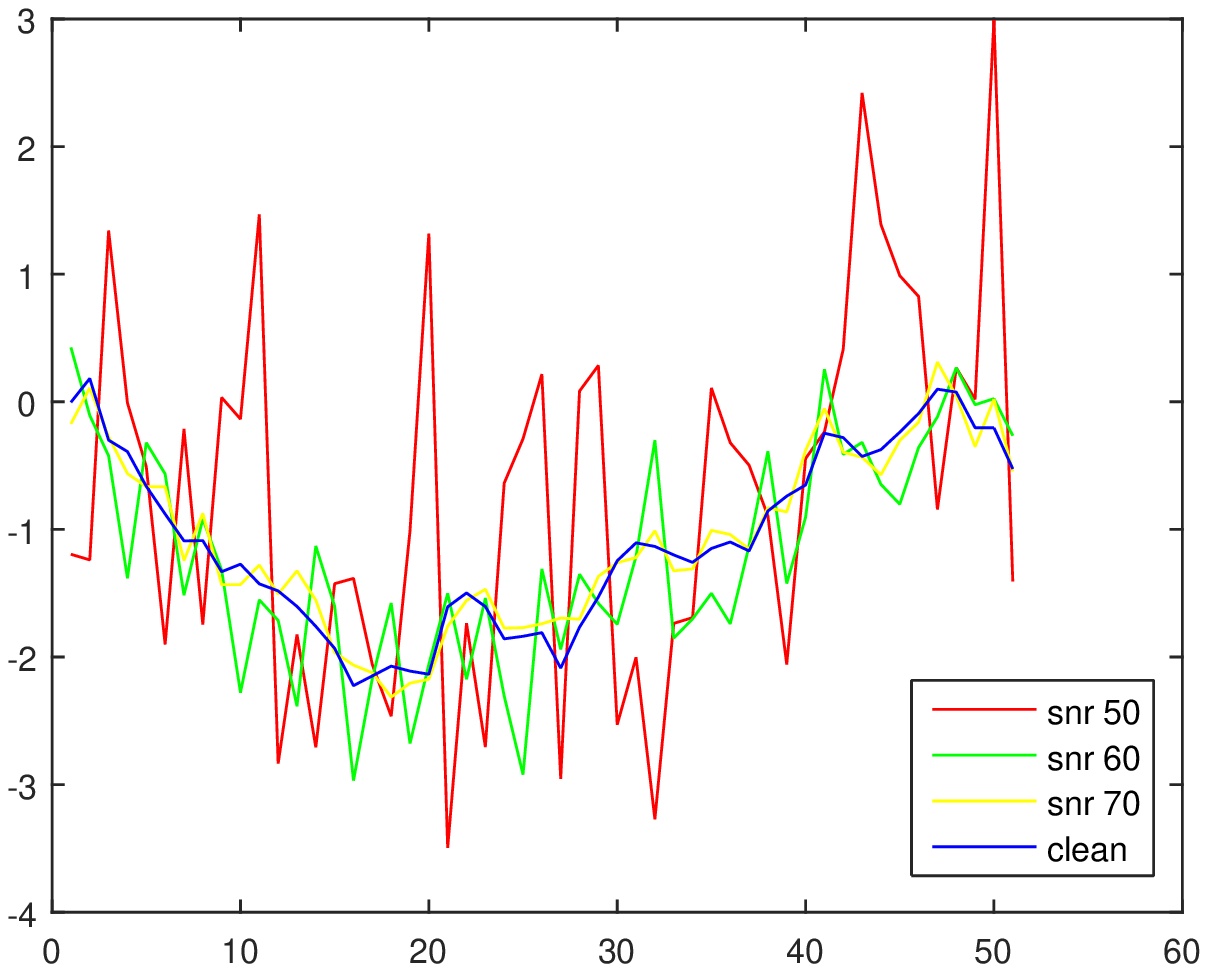}
\caption{Comparision of clean vs. noisy signals for 50 points of Wf ($D=1.2$)} 
\label{wnoise}
\end{figure}

\begin{figure}
\includegraphics[width=0.5\textwidth]{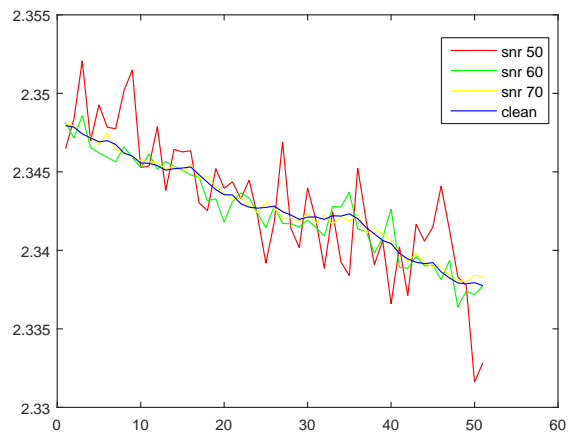}
\caption{Comparision of clean vs. noisy signals for 50 points of fBm ($D=1.2$)} 
\label{fnoise}
\end{figure}

In case of both Wf and fBm we are going to analyse the accuracy of the methods for different conditions:  

\begin{enumerate}
\item different length of the time--series (500, 10000,  or 314160 (Wf), 320000 (fBm) data points)
\item different sampling (full, each $10$th data point, each $100$th data point)
\item added white noise to data (using the function $awgn$ implemented in MatLab), with different signal-to-noise ratios per sample in dB: 50, 60, 70 (see Figure \ref{wnoise} and \ref{fnoise})
\end{enumerate}
In case of shorter non-overlapping parts of signal, $D$ was estimated for each part and the median of these was taken as the resulting estimate of $D$. 

\begin{figure}
\includegraphics[trim=50 20 50 10,clip,width=1\columnwidth,, height=6cm]{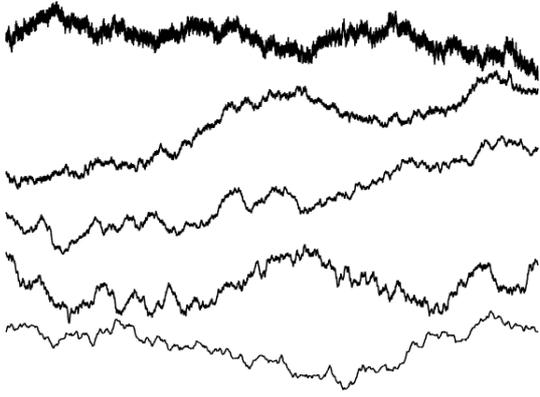}
\caption{Graphs of fBm with $D=1.2$, $D=1.4$, $D=1.5$, $D=1.6$, $D=1.8$ from bottom to top.} 
\label{diff} 
\end{figure} 
 
\section{Methods}

In this section the testing methods are described in detail. 

\subsection{Higuchi fractal dimension}

Higuchi’s approach of determining $D$ has been introduced in \cite{higuc1988}. The method scans fluctuations of the signal by investigating the defined length of the curve for different magnifications of the time axis of the signal. We take time-series of length $N: X(1),X(2),...,X(N)$ and make $k$ lagged time series that start from $m$-th place $(m=1,2,...,k)$ with gap of the size $k$:
$$X^{m}_k = X(m), X(m+k),...,X\left(m+\left[\frac{N-m}{k}\right]k\right)$$
Then Higuchi defines a length of the curve of $X^{m}_k$ as follows:
$$L^{m}_k =\sum_{i=1}^{\left[\frac{N-m)}{k}\right]}|X(m+ik)-X(m+(i-1)k|\frac{N-1}{\left[\frac{N-m}{k}\right]k^{2}}$$
When the length of the curve is calculated (and normalized) for every $m$ and $k$, we get an $L(k)$ as the mean of all lengths $L^{m}_k$. If $L(k)\sim k^{-D}$, the curve is a fractal with  dimension $D$. In case of real data application, we have to determine how high should the value $k$ go. In various articles, the recommendation regarding $k$ is to compute the estimates for increasing values of $k$ and use the value where the estimates reach a plateau. However, for our data we did not observe a plateau and in the end we made a heuristic choice of the value of $k=15$.

To estimate the Higuchi dimension, we used a Matlab code shared by J. Monge-Álvarez \cite{monge2014}.

\subsection{Generalized Hurst exponent} 
 
Originally, studies involving Hurst exponent were mainly connected to hydrology. H. E. Hurst (1880 – 1978) invented the exponent when determining the optimum dam sizing for the Nile river. Value of $H$ describes a dependence (or persistence) of the long-term memory process. The exponent reflects whether the time-series is persistent ($0.5<H<1$) or anti-persistent ($0<H<0.5$). That means how current data points influence future data points, whether the process is trending or mean reverting. 

As stated above, for self-affine time series $$D=2-H$$
Hurst\lq{s} estimation of $H$ was connected to so-called R/S statistics, where R is the range of partial sums of the time series and S is the sample standard deviation \cite{hurst1951}. 

In this study, the generalized Hurst exponent as described in \cite{dimat2003} was estimated. This approach begins with the q-order moments of the distribution of the increments: $$K_q (\tau)=\frac{<|X(t+\tau)-X(t)|^q>}{<|X(t)|^q>}$$

Then $H$ can be calculated from the following relation: $$ H(q)\sim \frac{\log K_q(\tau)}{q\log\tau}$$

$H(q) = H$ is characteristic of uni-scaling or uni-fractal processes, while $H(q)$ dependent on $q$ indicates processes called multi-scaling or multi-fractal. Different exponents $H(q)$ are associated with special features. For instance, $H(1)$ describes the scaling behaviour of the absolute values of the increments.
Depending on the context of the application, we have to set the value of $q$ and the range of $\tau$ values. For our purposes, $q=1$ and $ \tau$ from 1 to 25 were chosen.

\subsection{Detrended fluctuation analysis} 

DFA offers an alternative method of estimating the Hurst exponent \cite{peng1994}. 

The obtained exponent $\delta$ corresponds to the $H$. However, DFA has the advantage that it may also be applied to non-stationary signals.

DFA method looks on the fluctuations of the signal by investigating the variance of data segments on different scales. With rising scale the variance of data points decreases and the "speed" of the decrease determines $H$ and consequently $D$. 
To be more specific, we take  N data points: $ X(1),X(2),...,X(N)$ and make cumulative sum of the signal ($Y_{m}=\sum_{i=1}^{m}X(i)$). In the next step we divide the entire sequence into $N:l$ non-overlapping boxes of the size of $l$. Then we define the "local trend" $T(j)  (j=1,2,..., l)$ in each box to be the ordinate of a linear least-squares fit for the data points in box. In the next step every point is detrended by the belonging "local trend" and the variance for each box is calculated. The average of these variances over all boxes of size $l$ gives the value $F(l)$. For sufficiently high $l$, we observe the following relation: $$F(l)\sim l^\delta.$$ 
Finally, we make a log-log graph to find the slope  $\delta$. If $0<\delta<1$ the time-series is produced by a stationary process and $H=\delta$, if $1<\delta<2$ the process is non-stationary and $H=\delta -1$  \cite{peng1994}. The size of the boxes which we investigate is also important for accuracy of the result. We should not use too small sizes $(l<4)$, because finding the trend in small boxes would not be accurate enough and we also should avoid too big size $(l>N/4)$, because of insufficient number of  $F(l)$ values to average (less than 4) \cite{peng1994}. Since the sizes do not need to be divisors of $N$, in majority of cases there is some unused part at the end of signal. To take that part in consideration, for every $l$ we make two computation, one with the boxes starting from the beginning of the signal, and one starting from the end.

\subsection{Spectral decay coefficient} 
Discrete Fourier transform represents data by a superposition of sines and cosines that have various amplitudes and frequencies. With time series of length $N$, the range of frequencies that can be considered goes from $1/N$ to $1/2$ (Nyquist frequency given by the sampling). Powers of individual Fourier component are called power spectral density (power spectrum). 
They are usually computed by the technique called fast Fourier transform ( We used the \textit{fft} function from Matlab). For self-affine signals the power spectrum $P(f)$ decays via a power law \cite{mande1977}. 
$$P(f)\sim f^{-\beta} \to \beta=\frac{-\log P(f)}{\log f}$$
The value of the power-law decay can be obtained as a slope of linear regression applied to the power spectrum in log–log coordinate system.

After calculating the spectral decay coefficient $\beta$, we use the relationship introduced above to obtain the fractal dimension:
$$D=\frac{5-\beta}{2}$$

The power decay was estimated from the whole accessible frequency domain. However, use of restricted frequency range and also use of alternative spectrum estimators can be considered. All in all, successful use of this method requires a thorough knowledge of the pitfalls of the spectral analysis.

\begin{table*}[!h]
\includegraphics[trim=0 320 0 0,clip,width=1\textwidth]{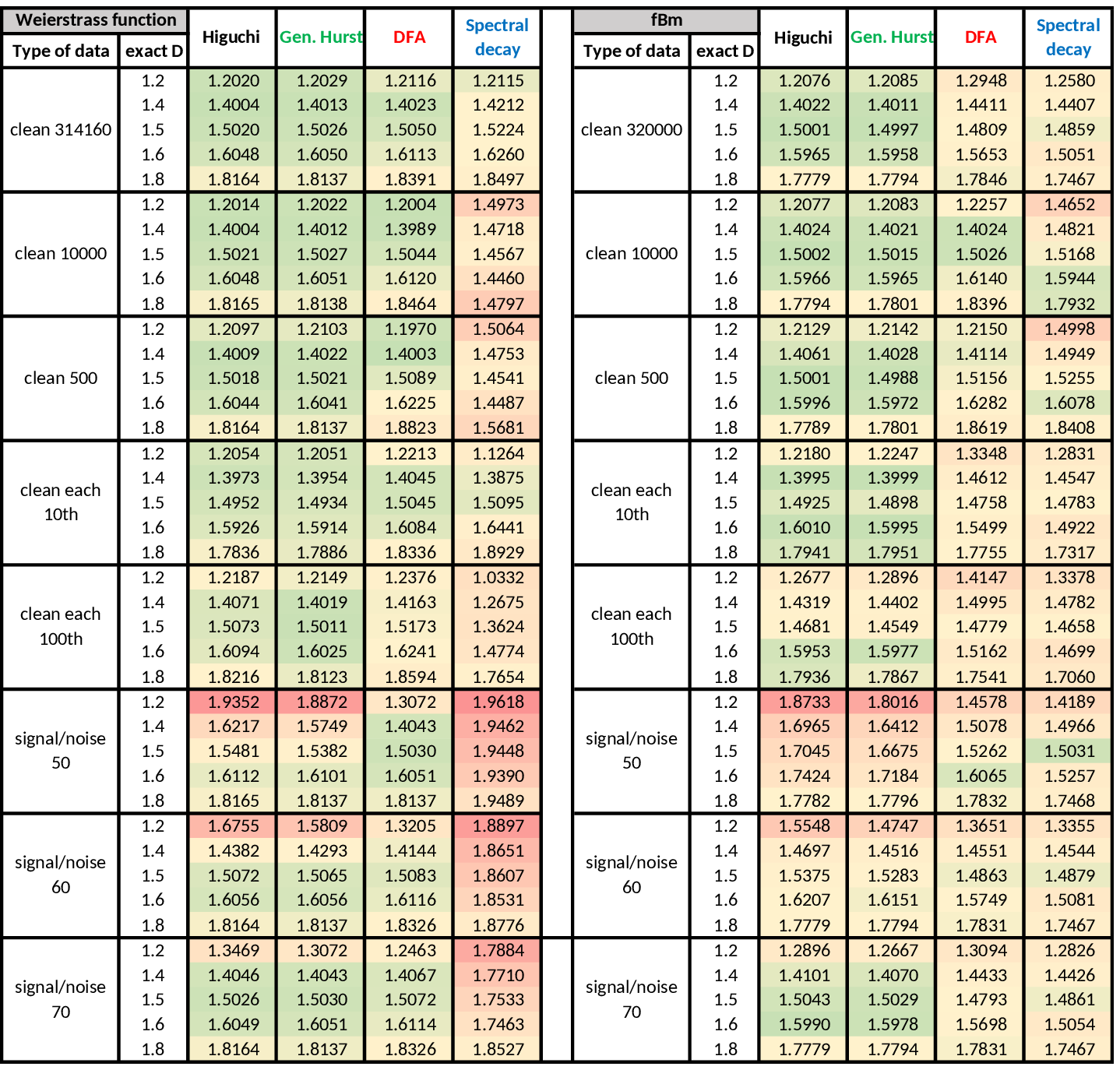}
\caption{Table of results for all types of tested data, methods and dimensions. The darker green, the more precise result; the darker red, the less precise result.} 
\label{tabulka}
\end{table*}  
\clearpage
\begin{figure}
\includegraphics[width=0.5\textwidth]{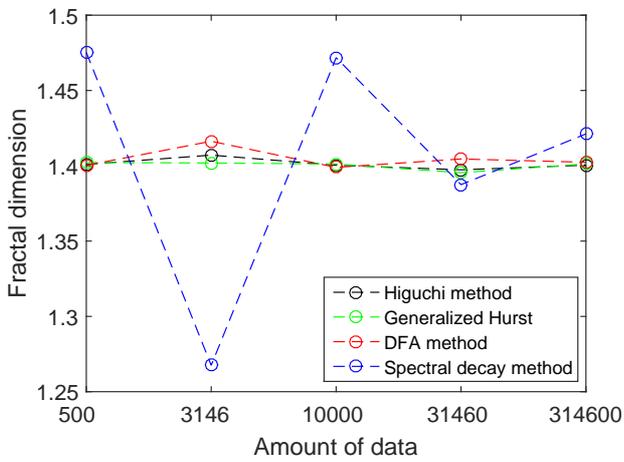}
\caption{Average results of estimation of $D$ depending on number of data for Weierstrass function. (Values 3141 and 31416 belong to coarse sampling - each 100th and 10th data points.)}
\label{ADw} 
\end{figure}  

\begin{figure}
\includegraphics[width=0.5\textwidth]{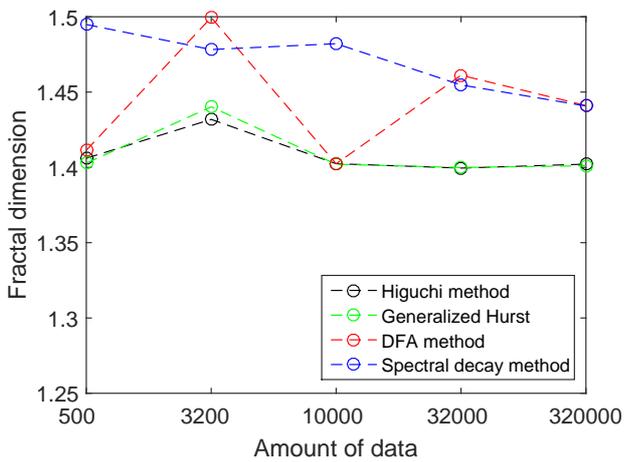}
\caption{Average results of estimation of $D$ depending on number of data for fBm.  (Values 3141 and 31416 belong to coarse sampling - each 100th and 10th data points.)}
\label{ADf} 
\end{figure}  
\begin{figure}
\includegraphics[width=0.5\textwidth]{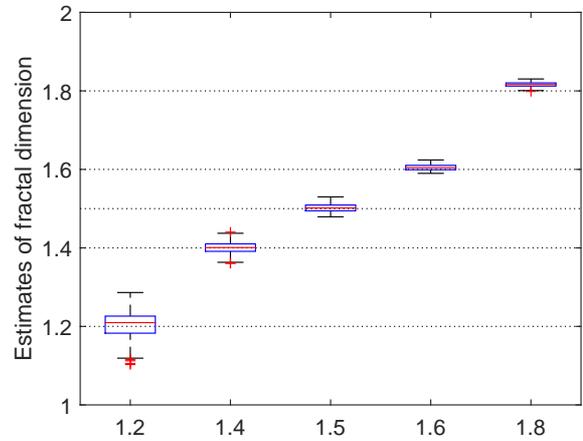}
\caption{Box plot for results of estimates of $D$, obtained by Higuchi method computed for $628$ segments of $500$ data points of Weierstrass function.}
\label{500w} 
\end{figure}  

\begin{figure}
\includegraphics[width=0.5\textwidth]{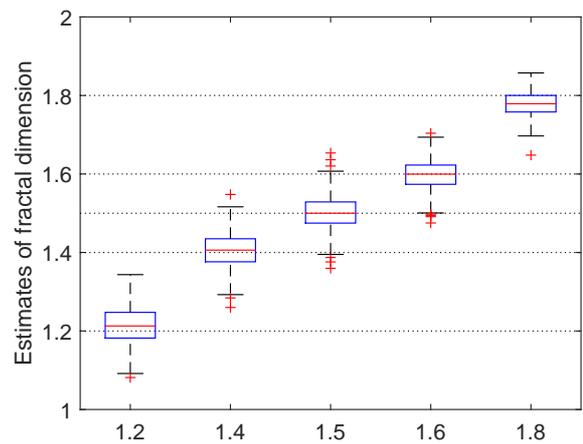}
\caption{Box plot for results of estimates of $D$, obtained by Higuchi method computed for $640$ segments of $500$ data points of fBm.}
\label{500f} 
\end{figure} 
\begin{figure}
\includegraphics[width=0.5\textwidth]{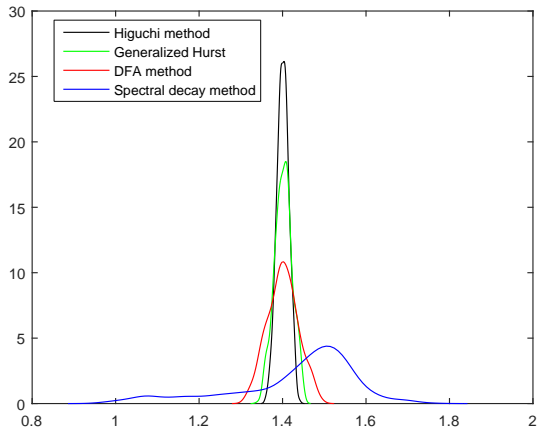}
\caption{Probability density of $628$ estimates of $D$, obtained by the four methods for $500$ point long segments of Weierstrass function ($D=1.4$).} 
\label{wdens}
\end{figure}  
\begin{figure}
\includegraphics[width=0.5\textwidth]{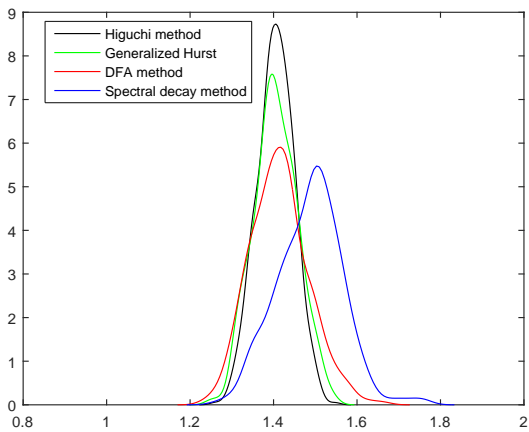}
\caption{Probability density of $640$ estimates of $D$, obtained by the four methods for $500$ point long segments of fBm ($D=1.4$).} 
\label{fdens}
\end{figure}  


\clearpage

\section{Results}
The results of this study are summarized in Table  \ref{tabulka}.

We can see that two of the implemented methods, the Higuchi method and the generalized Hurst exponent, performed generally the best. They gave comparably good $D$ estimations for both long and short data sets. Even for as short data sets as $500$ points and coarse sampling, the average results were very good (see Figures \ref {ADw}, \ref{ADf}, \ref{500w}, \ref{500f}, \ref{wdens}, \ref{fdens}). 

Detrended fluctuation analysis was not so successful as the best two methods. However, for noisy data, DFA unlike the other methods, seems to stay close to the dimension of the original signal, omitting the added complexity of the noise. 

The fourth method, the method of spectral decay, was the least successful. For this approach precise computation of the spectrum and optimum region for fitting the regression line are crucial but difficult to find. 

And finally, we can point out yet another problem: The production of a truly fractal signal of given dimension is itself a non-trivial task. So actually, to a certain extent we can worry about the adequacy of the data-generating algorithms as well as about the estimation algorithms.

\section{Conclusion}

Four methods for analysing self-affine fractal signals were tested in context of estimation of the fractal dimension.

If applied to our data, the Higuchi method and the method for generalized Hurst exponent were the best, while the DFA and especially the spectral decay gave more biased results. 

The results show that finding the optimal method is not easy. Performances of the four methods were not radically different.
Our recommendation could be the Higuchi method, partly because of the accuracy of the estimates and also because, unlike the other methods, it does not require preliminary adjustments of several sensitive parameters.

On the other hand, we do not rule
out possible improving of the performance of the remaining methods, especially after further investigations of how to optimize the parameter settings.

\section*{Acknowledgments}
This work was supported by the Slovak Grant Agency for Science (Grant no. 2/0011/16).

\bibliographystyle{apalike}
\bibliography{referencesm}

\end{document}